\documentclass[]{article}

\title{Lower bounds on the maximum dimension of a simple module in characteristic $p$}
\author{Geoffrey R. Robinson}

\usepackage{latexsym,amssymb}
\begin{document}

\maketitle

\begin{abstract}
	We obtain lower bounds for the maximum dimension of a simple $FG$-module, where $G$ is a finite group, $p$ is a prime, and $F$ is an algebraically closed field of characteristic $p$. These bounds are described in terms of properties of $p$-subgroups of $G$. It turns out that it suffices to treat the case that 
	$O_{p}(G) = \Phi(G) = 1.$ In that case, our main result asserts that if $p$ is an odd prime which is not Mersenne, then there is a simple $FG$ module $S$ of dimension at least $\frac{|G|_{p}}{|{\rm Out}_{G}(X)|_{p}}$,\\ where $X = O^{p^{\prime}}(E(G))$ and ${\rm Out}_{G}(X) = G/XC_{G}(X)$ is (isomorphic to) the group of outer automorphisms of $X$ induced by the conjugation action of $G$. Note, in particular, that this quotient is $|G|_{p}$ when $X = 1$, that is to say, when no component of $G$ has order divisible by $p$ (and, more generally, the quotient is $|G|_{p}$ when ${\rm Out}_{G}(X)$ is a $p^{\prime}$-group).
	
	\medskip
	We also use simplicial complex methods to obtain a weaker general bound which still applies when $p = 2$ or when $p$ is a Mersenne prime.

\end{abstract}

\medskip
\section{ Introduction and Statement of  Theorem 1}

\medskip
In this note, we are interested in exhibiting simple $FG$-module $S$ whose dimension is large relative to the size of $p$-subgroups of $G$, where $F$ is an algebraically closed field of characteristic $p$. We will prove:

\medskip
\noindent {\bf Theorem 1:} \emph{ Let $G$ be a finite group, $p$ be a prime, and $F$ be an algebraically closed field of characteristic $p$. Suppose further that $O_{p}(G) = \Phi(G) = 1$. Set $X = O^{p^{\prime}}(E(G))$ and ${\rm Out}_{G}(X) = G/XC_{G}(X).$ Then:}
	
\medskip	
\noindent \emph{i) 	If $p$ is an odd prime which is not Mersenne, then there is a simple $FG$ module $S$ of dimension at least $\frac{|G|_{p}}{|{\rm Out}_{G}(X)|_{p}}.$ } 

\medskip\noindent ii) \emph{ If $p=2$, or if $p$ is a Mersenne prime, then there is a simple $FG$-module of dimension at least $|A|$, where $A$ is a maximal Abelian $p$-subgroup of $XC_{G}(X)$ which contains an Abelian $p$-subgroup of maximal order of $C_{G}(X)$.}	

\medskip
\noindent {\bf Remark:}	Note, in particular, that the quotient given in part i) of Theorem 1 is $|G|_{p}$ when $X = 1$, that is to say, when no component of $G$ has order divisible by $p.$

\section{On the maximal dimension of a simple module}
 
 \medskip
 Let $G$ be a finite group, $p$ be a prime, and $F$ be an algebraically closed field of characteristic $p$. Let $P$ be a Sylow $p$-subgroup of $G$. We define the invariant $m_{s}(G) = m_{s}(G,p)$ to be the maximal dimension of a simple $FG$-module, and we are concerned in this note with obtaining lower bounds for $m_{s}(G)$ in terms of properties of $P$ and its subgroups. Clearly, we have $m_{s}(G) = m_{s}(G/O_{p}(G))$, so it suffices to deal with the case that $O_{p}(G) = 1.$ Note that we clearly have $m_{s}(G) = m_{s}(X)m_{s}(Y)$ when $G = X \times Y$ is a direct product of subgroups $X$ and $Y$.
 
 \medskip
 There are many examples of finite groups $G$ with $O_{p}(G) = 1$ and\\ $m_{s}(G) = |P|.$ For example, when $p = 2$ and $G =  {\rm SL}(2,2^{n}),$ then we have $m_{s} (G) = 2^{n}$,  the degrees of the simple $FG$-modules being $1,2,2^{2}, \ldots, 2^{n}$ (with $\left( \begin{array}{clcr} n\\ j \end{array}\right)$ occurences of simple modules of dimension $2^{j}$), so we can't expect to obtain a general bound stronger than $m_{s}(G) \geq |P|$. It is known that finite simple groups of Lie type have a $p$-blocks of defect zero for every prime $p$ by theorems of Michler [5] and Willems [9], so we do have $m_{s}(G) \geq |P|$ for such groups $G$ (whatever the defining characteristic of $G$ is). For most primes $p$ (other than the defining characteristic and prime divisors of the order of the Weyl group), such groups $G$ have Abelian Sylow $p$-subgroups. 
 
 \medskip
 It is clear that $m_{s}(G) \geq m_{s}(H)$ for each section $H$ of $G$. For if $N \lhd G$, we certainly have $m_{s}(G/N) \leq m_{s}(G)$, while if $K$ is a subgroup of $G$, and $S$ is a simple $FK$-module of maximal dimension, then ${\rm Ind}_{K}^{G}(S)$ has a composition factor $T$ such that $S$ is a composition factor of ${\rm Res}^{G}_{K}(T)$, so that\\ ${\rm dim}_{F}(T) \geq {\rm dim}_{F}(S).$
 
 \medskip
 Slightly less obvious than the reduction to $O_{p}(G) = 1$ is that it suffices to deal with the case $\Phi(G) = 1$ when $O_{p}(G) = 1.$ For if $O_{p}(G) = 1$, then $\Phi(G)$ is a nilpotent $p^{\prime}$-group. We reproduce the well-known proof that\\ $O_{p}(G/\Phi(G)) = 1.$ For suppose otherwise, and let $N \lhd G$ denote the full pre-image in $G$ of $O_{p}(G/\Phi(G)).$ Let $Q$ be a Sylow $p$-subgroup of $N$. Then $G = NN_{G}(Q)$ by the Frattini argument, and $N = \Phi(G)Q$, so that $G = \Phi(G)N_{G}(Q)$. Hence $G = N_{G}(Q)$, contrary to $O_{p}(G) = 1.$ 
 
 \medskip
 Note that when $\Phi(G) = 1$, we have $\Phi(M) = 1$ whenever $M \lhd \lhd  G $. In particular, each component of $G$ (if any exist) is a non-Abelian simple group.
 Also, $F(G)$ is a direct product of minimal normal subgroups of $G$, each of which is complemented in $G$. In particular, $F(G)$ is a direct product of elementary 
 Abelian groups (for various primes). Hence $F^{\ast}(G) = {\rm soc}(G)$ is the direct product of the minimal normal subgroups of $G$.
 
 \medskip
 In determining a lower bound for $m_{s}(G)$ in terms of properties of $p$-subgroups of $G$ when $O_{p}(G) = \Phi(G) = 1$, it suffices to treat the case that $G = PF^{\ast}(G)$.
For $O_{p}(PF^{\ast}(G))$ centralizes $F^{\ast}(G)$ in that case, and $Z(F(G))$ is a $p^{\prime}$-group. Hence we assume from now on that 
$O_{p}(G) = \Phi(G) = 1$, and that $G = PF^{\ast}(G),$ with $F^{\ast}(G)$ being a direct product of (possibly Abelian) simple groups.

\medskip
\section{ When $G \neq O^{p}(G)$}

\medskip
In the previous section, we reduced to a configuration where $G = PF^{\ast}(G)$. Note that $G \neq O^{p}(G)$ if $G \neq F^{\ast}(G).$ It may be worthwhile to digress to make some general remarks about the relationship between $m_{s}(G)$ and $m_{s}(H)$ when $H = O^{p}(G).$ In one sense, the Clifford theory between simple $FH$-modules and simple $FG$-modules is transparent, yet the relationship between $m_{s}(G)$ and $m_{s}(G)$ is less clear in general.

\medskip
Given a simple $FH$-module $S$, with inertial subgroup 
$I_{G}(S) = K$, say, then $S$ has a unique extension $\tilde{S}$ to a simple $FK$-module, and then ${\rm Ind}_{K}^{G}(S)$ is the unique simple $FG$-module covering $S$.
Hence the number of simple $FG$-modules is the number of orbits of $G$ on isomorphism types of simple $FH$-modules, and a simple $FH$-module $S$ in a given $G$-orbit gives rise to a simple $FG$-module of dimension $[G:I_{G}(S)]{\rm dim}_{F}(S)$.

\medskip
Hence we obtain $m_{s}(H) \leq m_{s}(G) \leq [G:H]m_{s}(H)$ but it usually difficult to ascertain without further knowledge which simple $FH$-module gives rise in this fashion to a simple $FG$-module of maximal dimension. 

\medskip
A case we will return to later is when $G = F^{\ast}(G)P$, and $G$ has no component of order divisible by $p$. Even in the case that $F^{\ast}(G) = N$ is Abelian, there are examples where $m_{s}(G) = |A|$, the maximal order of an Abelian $p$-subgroup of $G$.

\medskip
When $G = NP$, with $N$ an Abelian normal $p$-complement acted on faithfully by the finite $p$-group $P$, the discussion at the beginning of this section shows that each simple $FG$-module has dimension a power of $p$ dividing $|P|$. Also, each projective indecomposable $FG$-module has dimension $|P|$. In general, there need not be a projective simple $FG$-module (indeed, by dimension considerations, there certainly can be no such simple module if $|P| > |N|$, a situation which can occur).

\medskip
Examples of the last phenomenon occur as familiar types of counterexamples to ``Burnside's other $p^{a}q^{b}$-theorem", or to various types of regular orbit theorems.
One standard type is when $p$ is a Mersenne prime, $P \cong C_{p} \wr C_{p}$, and $N$ is an elementary Abelian group of order $(p+1)^{p}$. In this case, we do have $|P| > |N|$, so that $m_{s}(G) < |P|$, and we find that $m_{s}(G) = \frac{|P|}{p}$, which is the maximal order of an Abelian subgroup of $P$.

\medskip
Another standard type of counterexample is when $p = 2$ and $N$ is an elementary Abelian $q$-group of order $q^{2}$, where $q > 3$ is a Fermat prime.
Then the $2$-group $P = C_{q-1} \wr C_{2}$ acts faithfully on $N$, and we have $|P| > |N|$, the semidirect product $G = NP$ has $m_{s}(G) < |P|$, and we again 
find that $m_{s}(G) = \frac{|P|}{p}$, which is the maximum order of an Abelian subgroup of $P$. Even in the case $p=2,q = 3,$ we find that $m_{s}(G) = 4$, the maximum order of an Abelian $2$-subgroup of $G = NP$, although we do have $ |P| = 8 < 9 = |N|$ in this case. 

\medskip
When $p = 2$ or $p$ is a Mersenne prime, we may use direct products of some of the groups above to construct arbitrarily large examples of solvable $p$-nilpotent groups $G$ with $O_{p}(G) = 1$ such that $m_{s}(G)$ is the maximum order of an Abelian $p$-subgroup $A$ of $P$, yet $G$ has a non-Abelian Sylow $p$-subgroup $P$. In fact, it is clear that $[P:A]$ may be made as large as desired.

\medskip
In the next section, we use $p$-subgroup complexes to obtain a weaker general bound which suffices to cover these exceptional cases, but which will be significantly improved in the last section for odd primes $p$ which are not Mersenne.

\medskip
\section{Using $p$-subgroup complexes}
When $G$ is a finite group, we let $\mathcal{P}$ denote the simplicial complex associated to the poset of non-trivial $p$-subgroups of $G$. This was introduced by D. Quillen in [7], and has been extensively studied by K.S. Brown, S. Bouc, P.J. Webb, and J. Th\'evenaz, among others. We let $\sigma$ denote a strictly increasing chain of non-trivial $p$-subgroups of $G$. We let $|\sigma|$ denote the number of inclusions in the chain. The group $G$ acts by conjugation on chains, and the reduced Euler characteristic of this complex is (up to a sign depending on conventions chosen) $\sum_{\sigma \in \mathcal{P}/G} (-1)^{|\sigma|} [G:G_{\sigma}]$, where we include the empty chain, and consider it to have length $0$.

\medskip
In the Green Ring for $FG$, there is defined what P.J. Webb described as the Steinberg (virtual) module for $\mathcal{P}$, which is 
$$\sum_{\sigma \in \mathcal{P}/G} (-1)^{|\sigma|} {\rm Ind}_{G_{\sigma}}^{G}(F),$$ and which we denote by $St_{p}(G).$ 
This was proved in (Webb,[8]), to be a virtual projective module, that is to say, a difference of projective modules (allowing the possibility of the zero module).

\medskip
The virtual module ${\rm St}_{p}(G)$ may instead be calculated with respect to other $G$-homotopy equivalent complexes (with the same result). One of these, favoured by Quillen, is the complex $\mathcal{E}$ associated to the poset of elementary Abelian $p$-subgroups of $G$. Another, favoured by S. Bouc, is the complex $\mathcal{U}$ associated to the poset of non-trivial $p$-subgroups $U$ of $G$ which satisfy $U = O_{p}(N_{G}(U))$.

\medskip
Since ${\rm St}_{p}(G)$ is a virtual projective module, it is a difference of projective $FG$-modules, each of which is liftable to characteristic zero. Hence ${\rm St}_{p}(G)$ yields a well-defined complex virtual character. Furthermore, ${\rm St}_{p}(G)$ is uniquely determined (as a linear combination of indecomposable projectives) by this virtual character, using the non-singularity of the Cartan matrix. 

\medskip
If $O_{p}(G) \neq 1$, it is well-known (see, eg, Quillen [7]), that ${\rm St}_{p}(G) = 0$. We know of no example of a finite group with $O_{p}(G) = 1$ and ${\rm St}_{p}(G) = 0$, but it seems to be an open question at present whether $O_{p}(G) = 1$ implies ${\rm St}_{p}(G) \neq 0.$ 

\medskip
The virtual character afforded by (the lift to characteristic zero) of ${\rm St}_{p}(G)$  is the alternating sum of characters afforded by homology modules associated to the complex $\mathcal{P}$. In particular, if ${\rm St_{p}}(G) \neq 0,$ then the homology of the complex $\mathcal{P}$ is non-zero, and the complex $\mathcal{P}$ is not contractible.

\medskip
In cases where we are able to  establish that ${\rm St}_{p}(G) \neq 0$, we obtain a lower bound for $m_{s}(G):$ 

\medskip
\noindent {\bf Lemma 1:} \emph{Let $G$ be a finite group such that ${\rm St}_{p}(G) \neq 0$. Then there is simple $FG$-module $S$ and a chain $\sigma$ of non-trivial $p$-subgroups of $G$, such that ${\rm Res}^{G}_{G_{\sigma}}(S)$ contains the projective cover of the trivial module as a summand. In particular,
	${\rm dim}_{F}(S) \geq |Q|$, where $Q \in {\rm Syl}_{p}(G_{\sigma}).$ The chain $\sigma$ may be chosen to consist of elementary Abelian $p$-subgroups of $G$, in which case we have $|Q| \geq |A|$ for some maximal Abelian $p$-subgroup $A$ of $G$.}

\medskip
\noindent {\bf Proof:} Let $S$ be a simple module such that the projective cover of $S$ occurs with non-zero multiplicity in ${\rm St}_{p}(G).$ Then there is a chain $\sigma$ such that the projective cover of $S$ occurs as a summand of ${\rm Ind}_{G_{\sigma}}^{G}(F),$ and $\sigma$ may be chosen to consist of elementary Abelian $p$-subgroups.

\medskip
Then the projective cover of the trivial module occurs as a summand of ${\rm Res}^{G}_{G_{\sigma}}(S)$, since ${\rm Ind}_{G_{\sigma}}^{G}(F \otimes {\rm Res}^{G}_{G_{\sigma}}(S^{\ast}))$ has the projective cover of the trivial module as a direct summand, and the multiplicity of the projective cover of the trivial module as a summand is preserved by induction of modules.

\medskip
Hence the first claim follows. If we choose $\sigma$ to consist of elementary Abelian subgroups, and we choose a maximal Abelian $p$-subgroup $A$ of $G$ containing the largest subgroup of the chain, then $A \leq G_{\sigma}$ and $|Q| \geq |A|$.

\medskip
Next, we observe that when $F^{\ast}(G)$ is a $p^{\prime}$-group, (in particular, when $G$ is $p$-solvable with $O_{p}(G)= 1$), then we can do somewhat better.

\medskip
\noindent {\bf Corollary 2:} \emph{ Let $G$ be a finite group such that $F^{\ast}(G)$ is a $p^{\prime}$-group, and let $A$ be an Abelian $p$-subgroup of $G$ of maximal order. Then $m_{s}(G) \geq |A|$.}

\medskip
\noindent {\bf Proof:} It suffices to prove that $m_{s}(O_{p^{\prime}}(G)A) \geq |A|$, so we may suppose that $G = O_{p^{\prime}}(G)A$, since 
$O_{p}(O_{p^{\prime}}(G)A) = 1.$

\medskip
But in this case, a Theorem of Hawkes and Isaacs ([4]) applies, and allows us to conclude that the virtual module ${\rm St}_{p}(G) \neq 0.$ In fact, Hawkes and Isaacs even prove that the reduced Euler characteristic of $\mathcal{P}$ is non-zero, and the reduced Euler characteristic of $\mathcal{P}$ is the virtual dimension of ${\rm St}_{p}(G).$ Hence Lemma 1 applies to $G$ (note that, under current assumptions, every maximal Abelian $p$-subgroup of $G$ is conjugate to $A$).

\medskip
In the paper [1] of Aschbacher-Kleidman, they establish a condition which, using a previous result of the present author, ensures that for any finite almost simple group $G$ with  $F^{\ast}(G) \not \cong {\rm PSL}(3,4)$, we have ${\rm St}_{p}(G) \neq 0$. They obtain the same conclusion
when $G$ is the simple group ${\rm PSL}_{3}(4)$, but the condition may fail for certain subgroups of ${\rm Aut}({\rm PSL}(3,4)).$ This yields:

\medskip
\noindent {\bf Corollary 3:} \emph{Let $G$ be a finite non-Abelian simple group. Then $m_{s}(G) \geq |A|$ for some maximal Abelian $p$-subgroup $A$ of $G$.}

\medskip
\noindent {\bf Proof:} By a Theorem of Aschbacher-Kleidman ([1]), we have ${\rm St}_{p}(G) \neq 0$.
Hence Lemma 1 may be applied to $G$, and $m_{s}(G) \geq |A|$ for some maximal Abelian $p$-subgroup of $G$.

\medskip
\section{ Proof of Theorem 1}

\medskip
By the results of Granville and Ono [2], Michler [5], and Willems [9], whenever $G$ is a finite non-Abelian simple group and $p$ is a prime greater than $3$, the group $G$ has a $p$-block of defect zero. However, for both $p = 2$ and $p = 3$, there are alternating groups and sporadic simple groups which have no $p$-block of defect zero.

\medskip
By the proof of Lemma 2.3 of Guralnick and Robinson [3], whenever $p$ is an odd prime which is not Mersenne, and $q \neq p$ is a prime, then a Sylow $p$-subgroup $R$ of 
$G = {\rm GL}(n,q)$ has at least two  regular orbits on the natural module for $G$. When $q=2$ and $p$ is not Mersenne, we need to adapt the latter proof as follows: In the case that $R$ is Abelian of order $p^{r}$ and the action of $R$ on the natural module for $G$ is irreducible, all orbits of $R$ on non-zero vectors are regular, but we don't have $2^{n}-1 = p^{r},$ since that equality would force $r$ to be odd and $p+1$ to be a power of $2$, contrary to the fact that $p$ is not a Mersenne prime. Hence there are at least two regular orbits of $R$ on the non-zero vectors of the natural module for $G$). Then the inductive argument proceeds in a manner similar to that of Lemma 2.3 of [3].

\medskip
The proof of our main result subdivides naturally according to whether $F^{\ast}(G)$ is a $p^{\prime}$-group or not.

\medskip
\noindent {\bf Lemma 2 :} \emph{ Let $G$ be a finite group such that $F^{\ast}(G)$ is a $p^{\prime}$-group, where $p$ is an odd prime which is not Mersenne. Then $m_{s}(G) \geq |P|$, where $P \in {\rm Syl}_{p}(G)$.}

\medskip
\noindent {\bf Proof:} As in earlier arguments, we may assume that $\Phi(G) = 1$ and that $G = PF^{\ast}(G)$. We may also assume that $P$ does not act faithfully on any proper $P$-invariant subgroup of $F^{\ast}(G)$. By Theorem 1.2 of Moret\'o-Navarro [6], we see that $F^{\ast}(G)$ is nilpotent (and even Abelian of squarefree exponent, since $\Phi(G) = 1$).

\medskip
Now it follows that $P$ has a regular orbit on the Abelian group\\ $F^{\ast}(G) = F(G)$, and on the dual group, which is the group of irreducible characters of $F(G)$. Hence $PF(G)$ has a $p$-block of defect zero and there is a simple $FG$-module of dimension $|P|$. Thus $m_{s}(G) \geq |P|$ (in fact, equality holds in this residual configuration).

\medskip
\noindent {\bf Conclusion of Proof of Theorem 1:} i) Suppose that $p$ is neither $2$ nor a Mersenne prime, and that $O_{p}(G) = \Phi(G) = 1.$ Let $X = O^{p^{\prime}}(E(G))$. Then $X$  is the normal subgroup of $E(G)$ generated the $p$-element of $E(G)$, and is the direct product of the components of $G$ of order divisible by $p$ (or is the identity subgroup if there is no such component) 

\medskip
It follows from remarks at the beginning of this section that $X$ has a $p$-block of defect zero, since $p \geq 5.$ Now $C \cap X = 1$, where $C = C_{G}(X)$, and we note that $F^{\ast}(C)$ is a $p^{\prime}$-group. By the previous Lemma, we have $m_{s}(C) \geq  |C|_{p}$ so that we have $m_{s}(G) \geq m_{s}(XC) \geq |XC|_{p}.$ Note that $C = G$ if $X = 1.$

\medskip
Now $G/XC$ is isomorphic to a subgroup of ${\rm Out}(X),$ and we denote this subgroup by ${\rm Out}_{G}(X).$ Now we have $m_{s}(G) \geq |XC|_{p} = \frac{|G|_{p}}{|Out_{G}(X)|_{p}}$.

\medskip
\noindent ii) If $p =2$ or a Mersenne prime, then we may argue in a fashion similar to i), but making use instead of Corollaries 2 and 3.

\medskip
\begin{center}
\bf{References}
\end{center}

\medskip
\noindent [1] Aschbacher, M.; Kleidman, P. : \emph{On a conjecture of Quillen and a Lemma of Robinson}, Arch. Math. (Basel) 55,(1990),3, 209-217.

\medskip
\noindent [2] Granville, A. ; Ono, K. : \emph{Defect zero blocks for finite simple groups }, Trans. Amer. Math. Soc., 348,1,(1996),331-347.

\medskip
\noindent [3] Guralnick, R. ; Robinson, G. : \emph{Some variants of the Brauer-Fowler Theorems}, J. of Algebra, 558,(2020),453-484.

\medskip
\noindent [4] Hawkes, T.; Isaacs, I.M.: \emph{ On the poset of $p$-subgroups of a $p$-solvable group}, J.London Math. Soc., 38,(1988), 77-86.

\medskip
\noindent [5] Michler, G. : \emph{ A finite simple group of Lie type has $p$-blocks with different defects, $p \neq 2$}, J. Algebra 104, (1986), 220-230.

\medskip
\noindent [6] Moret\'o, A. ; Navarro, G. : \emph{Heights of characters in blocks of $p$-solvable groups}, Bull. LMS, 37,(2008),373-380.

\medskip
\noindent [7] Quillen, D. : \emph{Homotopy properties of the poset of non-trivial $p$-subgroups},  Adv. in Math. 88, (1978),101-128.

\medskip
\noindent [8] Webb, P.J. : \emph{Complexes, group cohomology and an induction theorem for the Green ring}, J. Algebra 104,(1986), 351-357.

\medskip
\noindent [9] Willems, W. : \emph{Blocks of defect zero for finite simple groups of Lie type}, J. Algebra, 113, 2, (1988), 511-522.

\end{document}